\pdfoutput=1
\documentclass[a4paper]{scrartcl}
\usepackage{a4}
\usepackage[francais]{babel}
\usepackage[latin1]{inputenc}
\usepackage{amsfonts}
\usepackage{amssymb}
\usepackage{amsmath}
\usepackage[mathscr]{euscript}
\usepackage{theorem}
\usepackage[all]{xy}
\usepackage{vmargin}
\setmarginsrb{3cm}{2.5cm}{2.5cm}{3.5cm}{1cm}{1cm}{0cm}{0.8cm}

\renewcommand\leq{\leqslant}
\renewcommand\geq{\geqslant}

\newcommand\Nat{\mathbb{N}}

\newcommand\Z{\mathbb{Z}}
\newcommand\Q{\mathbb {Q}}
\newcommand\Fi{\mathbb{F}}

\newcommand\Hp{{}^{p}{\mathrm{H}}}

\newcommand\Aff{\mathbb{A}}

\newcommand\Gr{\mathbb{G}}
\newcommand\Bo{\mathcal{B}}

\newenvironment{preuve}{\vspace{0.3cm}\emph{Démonstration.}}{\hfill $\square$ \vspace{0.3cm}}
\newenvironment{preuvet}[1]{\vspace{0.3cm}\emph{Démonstration du théorème #1.}}{\hfill $\square$ \vspace{0.3cm}}

\DeclareMathOperator{\D}{D}
\newcommand\HH{\mathcal{H}}
\newcommand\Tate{\mathcal{T}}
\newcommand\DP{{}^{w}\D}

\newcommand\as{\underline{a}}
\newcommand\fl{\longrightarrow}

\newcommand\iso{\stackrel {\sim}{\fl}}

\theoremstyle{plain}

\newtheorem{theoreme}[subsection]{Théorème}
\newtheorem{prop}[subsection]{Proposition}
\newtheorem{proposition}[subsection]{Proposition}
\newtheorem{lemme}[subsection]{Lemme}
\newtheorem{souslemme}[subsubsection]{Sous-lemme}

\title{Note sur les polynômes de Kazhdan-Lusztig}
\author{Sophie Morel}

\begin{document}

\maketitle

Dans son article \cite{B}, Brenti a prouvé de manière combinatoire une
formule non-récursive pour les polynômes de Kazhdan-Lusztig d'un groupe
de Coxeter (cf le théorème 4.1 de \cite{B}).
Le but de cette note est de donner une interprétation géométrique 
de cette formule
pour les groupes de Coxeter qui sont isomorphes à un groupe de Weyl. Cette
interprétation utilise un résultat de \cite{M} (théorème 3.3.5), qui exprime le
prolongement intermédiaire d'un faisceau pervers pur comme ``tronqué par
le poids'' de son image.

Donnons rapidement le plan de l'article. La section 1 rappelle la
définition, due à Kazhdan et Lusztig, des polynômes de Kazhdan-Lusztig et
des polynômes $R$. La section 2 rappelle le résultat de \cite{M} qui sera
utilisé dans l'interprétation géométrique de la formule de Brenti.
La section 3 rappelle une décomposition des cellules de Bruhat due à
Deodhar, et en donne une nouvelle preuve, qui utilise des calculs de
Härterich. Cette décomposition sera utile dans la section 4, pour montrer
que la cohomologie des cellules de Bruhat a une forme particulièrement
simple. Enfin, la section 4 contient le résultat principal (théorème
\ref{calcul_P_R}), c'est-à-dire la preuve géométrique de la formule de
Brenti.

Je remercie Gérard Laumon de m'avoir signalé cette application de \cite{M}.

\vspace{1cm}

Dans toute la suite, $\Fi_q$ est un corps fini, $\Fi_q\subset\overline{\Fi}_q$ est une 
clôture algébrique de $\Fi_q$ et $\ell$ est un nombre premier inversible dans $\Fi_q$.

\section{Polynômes de Kazhdan-Lusztig et polynômes $R$}

Dans cette section, nous rappelons quelques résultats de Kazhdan et
Lusztig (cf \cite{KL1} et \cite{KL2}).

Soit $G$ un groupe algébrique réductif connexe déployé sur $\Fi_q$. Soient $B$
un sous-groupe de Borel de $G$ (défini sur $\Fi_q$), $T\subset B$ un tore maximal
déployé sur $\Fi_q$, $B^*$ le sous-groupe de Borel de $G$ opposé à $B$ (ie
tel que $B\cap B^*=T$),
$N$ le normalisateur de $T$ dans $G$ et $W=N/T$ le groupe de Weyl.
On note $\Phi$ l'ensemble des racines de $T$ dans $Lie(G)$, $\Phi^+$ le sous-ensemble de
racines positives associé à $B$ (c'est-à-dire l'ensemble des racines de $T$ dans
$Lie(B)$), et $\Delta\subset\Phi^+$ l'ensemble des racines simples. Enfin, on note $\ell$
la fonction longueur sur $W$ associée à $\Delta$ et $\leq$ l'ordre de Bruhat sur $W$.

Rappelons que l'algèbre de Hecke $\HH$ de $W$ est l'algèbre sur l'anneau 
$\Z[t^{1/2},t^{-1/2}]$ engendrée par des éléments $T_w$, $w\in W$, soumis aux relations :
$$\begin{array}{ll}T_wT_{w'}=T_{ww'} & \mbox{ si }\ell(ww')=\ell(w)+\ell(w') \\
(T_s+1)(T_s-t)=0 & \mbox{ si }s\mbox{ est la réflexion associée à une racine simple.}\end{array}$$

On définit des polynômes $R_{v,w}\in\Z[t]$, $v,w\in W$, $v\leq w$, par les formules
(cf \cite{KL1} $\S$ 2) :
$$T_w^{-1}=\sum_{v\leq w}(-1)^{\ell(w)-\ell(v)}R_{v,w}(t)t^{-\ell(w)}T_v.$$
$R_{v,w}$ est de degré $\ell(w)-\ell(v)$.

D'après le théorème 1.1 de \cite{KL1} (cf aussi \cite{KL2} $\S$ 2), 
il existe une unique famille de polynômes $P_{v,w}\in\Z[t]$,
$v,w\in W$, $v\leq w$, avec $P_{w,w}=1$ et $P_{v,w}$ de degré 
$\leq \frac{1}{2}(\ell(w)-\ell(v)-1)$ si $v<w$,
telle que : pour tous $v,w\in W$ tels que $v\leq w$,
$$t^{\ell(w)-\ell(v)}P_{v,w}(t^{-1})=\sum_{v\leq y\leq w}R_{v,y}(t)P_{y,w}(t).$$

\vspace{.5cm}

Rappelons l'interprétation géométrique des polynômes $P_{v,w}$ donnée par Kazhdan et
Lusztig dans \cite{KL2}.

Le choix de $B$ détermine un isomorphisme entre la variété $\Bo$ des sous-groupes de
Borel de $G$ et le quotient $G/B$. On a deux stratifications de $\Bo$ par des sous-variétés
localement fermées, $\Bo=\bigcup\limits_{w\in W}X_w=\bigcup\limits_{w\in W}X^w$, avec, pour
tout $w\in W$,
$$X_w=\{B'\in\Bo|\exists g\in B wB,B'=gB g^{-1}\}$$
$$X^w=\{B'\in\Bo|\exists g\in B^*wB,B'=gB g^{-1}\}.$$
La variété $X_w$ est isomorphe à l'espace affine $\Aff^{\ell(w)}$, et $X^w$ est isomorphe
à $\Aff^{dim(\Bo)-\ell(w)}$. Pour tout $w\in W$, on note $x_w=wBw^{-1}\in X_w\cap X^w$.

L'adhérence $\overline{X}_w$ de $X_w$ dans $\Bo$ est la variété de Schubert associée à
$w$. C'est une variété projective irréductible de dimension $\ell(w)$, et on a
$$\overline{X}_w=\bigcup_{v\leq w}X_v.$$
De même, on a
$$\overline{X}^w=\bigcup_{v\geq w}X^v.$$

On note $F:\Bo\fl\Bo$ le Frobenius. Pour tout $w\in W$, on a $F(X_w)\subset F(X_w)$, 
$F(X^w)\subset F(X^w)$ et $F(x_w)=x_w$.
Pour tous $v,w\in W$ tels que $v\leq w$, on note $j_w$ et $i_{v,w}$ les inclusions
de $X_w$ et $X_v$ dans $\overline{X}_w$, et
$$IC_{\overline{X}_w}=(j_{w!*}(\Q_\ell[\ell(w)]))[-\ell(w)]$$
le complexe d'intersection à coefficients constants de $\overline{X}_w$.

Kazhdan et Lusztig ont démontré le théorème suivant (\cite{KL2}, théorèmes 4.2 et 4.3) :

\begin{theoreme}\label{calcul_P_IC}(Kazhdan-Lusztig)
Soit $w\in W$. Alors le faisceau de cohomologie
(ordinaire) $\HH^i(IC_{\overline{X}_w})$ est nul si $i$ est impair. Si $i$ est pair et
$B'\in\overline{X}_w$ est stable par une puissance $F^r$ de $F$,
alors  les valeurs 
propres de $(F^r)^*$ sur la fibre
$\HH^i(IC_{\overline{X}_w})_{B'}$ sont toutes égales à $q^{ir/2}$.

De plus, pour tout $v\leq w$, on a
$$P_{v,w}(t)=\sum_{i\geq 0} dim(\HH^{2i}(IC_{\overline{X}_w})_{x_v})t^i.$$

\end{theoreme}

\section{Rappel d'un résultat de \cite{M}}

Soit $X$ un schéma séparé de type fini sur $\Fi_q$. On note $D^b_m(X,\Q_\ell)$ la catégorie
des complexes $\ell$-adiques mixtes sur $\Fi_q$. Cette catégorie est munie de la t-structure
donnée par la perversité autoduale, et on note $\Hp^i$ les foncteurs de cohomologie pour
cette t-structure. (Les notions de ce paragraphe sont définies dans
le chapitre 5 du livre \cite{BBD} de Beilinson, Bernstein et Deligne.)

Pour tout $a\in\Z\cup\{\pm\infty\}$, on note $\DP^{\leq a}(X)$ (resp. $\DP^{\geq a}(X)$) la
sous-catégorie pleine de $D^b_m(X,\Q_\ell)$ dont les objets sont les complexes mixtes $K$
tels que pour tout $i\in\Z$, le faisceau pervers $\Hp^iK$ soit de poids $\leq a$
(resp. $\geq a$). Alors (\cite{M}, proposition 3.1.1 (iv)) :

\begin{proposition} Pour tout $a\in\Z\cup\{\pm\infty\}$, $(\DP^{\leq a}(X),\DP^{\geq a+1}(X))$
est une t-structure sur $D^b_m(X,\Q_\ell)$.

\end{proposition}

On note $w_{\leq a}$ et $w_{\geq a+1}=w_{>a}$ les foncteurs de troncature pour cette 
t-structure. Ce sont des foncteurs triangulés (car $\DP^{\leq a}(X)$ et 
$\DP^{\geq a+1}(X)$ sont des sous-catégories triangulées de $D^b_m(X,\Q_\ell)$), 
et la dualité de Poincaré échange $w_{\leq a}$ et $w_{\geq -a}$.

\vspace{.5cm}

Soit $(S_0,\dots,S_n)$ une partition de $X$ par des sous-schémas (localement fermés non vides)
telle que, pour tout $k\in\{0,\dots,n\}$, $S_k$ soit ouvert dans 
$X-\bigcup\limits_{l<k}S_l$. En particulier, $S_0$ est ouvert dans $X$.
Pour tout $k\in\{0,\dots,n\}$, on note $i_k$ l'inclusion de $S_k$ dans $X$. 

Le théorème suivant découle des résultats de la section 3 de \cite{M} :

\begin{theoreme}\label{calcul_IC}
Soient $a\in\Z$ et $K$ un faisceau pervers pur de poids $a$ sur $S_0$.
Alors on a une égalité de classes
dans le groupe de Grothendieck de $D^b_m(X,\Q_\ell)$ :
$$\begin{array}{rcl}\displaystyle{[i_{0!*}K]} & = & \displaystyle{\sum_{1\leq n_1<\dots<n_r<n}(-1)^r[i_{n_r!}w_{\leq a}i_{n_r}^!\dots i_{n_1!}w_{\leq a}i_{n_1}^!i_{0!}K]} \\
 & & \displaystyle{+\sum_{1\leq n_1<\dots<n_r=n}(-1)^r[i_{n_r!}w_{<a}i_{n_r}^!i_{n_{r-1}!}w_{\leq a}i_{n_{r-1}}^!\dots i_{n_1!}w_{\leq a}i_{n_1}^!i_{0!}K].}\end{array}$$

\end{theoreme}

\begin{preuve} On utilise les notations de la section 3 de \cite{M}. D'après la proposition
\cite{M} 3.4.2, on a des isomorphismes canoniques 
$$w_{\geq (a,a+1,\dots,a+1)}i_{0!}K=w_{\geq (a,a,\dots,a)}i_{0!}K=i_{0!*}K,$$
d'où un isomorphisme canonique
$$w_{\geq (a,a+1,\dots,a+1,a)}i_{0!}K=i_{0!*}K.$$
De plus, le théorème \cite{M} 3.3.5 donne pour tout 
$\as=(a_0,\dots,a_n)\in (\Z\cup\{\pm\infty\})^{n+1}$ tel que $a=a_0$
une égalité dans le groupe de Grothendieck de $D^b_m(X,\Q_\ell)$ :
$$[w_{\leq\as}Ri_{0*}K]=\sum_{1\leq n_1<\dots <n_r\leq n}(-1)^r[Ri_{n_r*}w_{>a_{n_r}}i_{n_r}^*
\dots Ri_{n_1*}w_{>a_{n_1}}i_{n_1}^*Ri_{0*}K],$$
d'où par dualité :
$$[w_{\geq\as}i_{0!}K]=\sum_{1\leq n_1<\dots<n_r\leq n}(-1)^r[i_{n_r!}w_{<a_{n_r}}i_{n_r}^!
\dots i_{n_1!}w_{<a_{n_1}}i_{n_1}^!i_{0!}K].$$
Cette égalité, combinée avec l'isomorphisme ci-dessus, donne le résultat cherché.

\end{preuve}

\section{Décomposition des $X_w$}

Dans \cite{D2}, Deodhar a construit une décomposition des cellules de Bruhat $X_w$ qui induit
une décomposition agréable des $X_w\cap X^v$, pour $v\leq w$. Nous allons donner ici
une autre interprétation de cette décomposition.

Soit $w\in W$. On fixe une décomposition $w=s_1\dots s_r$ de $w$ en réflexions simples,
avec $r=\ell(w)$. Pour tout $i\in\{1,\dots,r\}$, on note $\alpha_i$ la racine simple
correspondant à $s_i$.
Soit $\Gamma=\{1,s_1\}\times\dots\times\{1,s_r\}$.
Pour tout $\gamma=(\gamma_1,\dots,\gamma_r)\in \Gamma$, on note 
$$I(\gamma)=\{i\in\{1,\dots,r\}|\gamma_i=s_i\}$$
$$J(\gamma)=\{i\in\{1,\dots,r\}|\gamma_1\gamma_2\dots\gamma_i(-\alpha_i)\in\Phi^+\}.$$
Pour tout $v\leq w$, on note
$$\Gamma_v=\{(\gamma_1,\dots,\gamma_r)\in\Gamma|\gamma_1\dots\gamma_r=v\}.$$

Deodhar a montré le résultat suivant (\cite{D2}, théorème 1.1 et corollaire 1.2) :

\begin{prop}\label{structure_link}(Deodhar)
La cellule de Bruhat $X_w$ s'écrit de manière
canonique comme une union disjointe de sous-variétés localement fermées
$$X_w=\bigcup_{\gamma\in\Gamma}Y_\gamma,$$
avec $Y_\gamma=\varnothing$ si $J(\gamma)\not\subset I(\gamma)$,
et, pour tout $\gamma\in\Gamma$ tel que $J(\gamma)\subset I(\gamma)$, 
$$Y_\gamma\simeq\Aff^{card(I(\gamma)-J(\gamma))}\times\Gr_m^{r-card(I(\gamma))}.$$
De plus, pour tout $v\in W$ tel que $v\leq w$, on a
$$X_w\cap X^v=\bigcup_{\gamma\in\Gamma_v}Y_\gamma.$$

\end{prop}

En utilisant des calculs de Härterich (\cite{H}, $\S$ 1), on peut donner une autre preuve de ce 
résultat :

Soit $\lambda:\Gr_m\fl T$ un cocaractère tel que $\langle\alpha,\lambda\rangle
>0$ pour tout
$\alpha\in\Phi^+$. On fait agir $\Gr_m$ sur $\Bo$ via $\lambda$. Les points fixes de cette
action sont les $x_z$, $z\in W$, et les deux décompositions de Bialynicki-Birula
(cf \cite{BB}) de $\Bo$ sont :
\begin{itemize}
\item[-] la décomposition en cellules contractantes, $\Bo=\bigcup\limits_{z\in W}X_z$;
\item[-] la décomposition en cellules dilatantes, $\Bo=\bigcup\limits_{z\in W}X^z$.

\end{itemize}

On va décomposer $X_w$ en utilisant la décomposition de Bialynicki-Birula de
la résolution de Bott-Samelson de $\overline{X}_w$ associée à la décomposition 
$w=s_1\dots s_r$.

Pour toute racine $\alpha$, on note $p_\alpha:\Aff^1\fl U_\alpha$ le sous-groupe
à un paramètre associé à $\alpha$.
Si $\alpha$ est une racine simple, on note $P_\alpha$ le sous-groupe parabolique (minimal)
de $G$ engendré par $B$ et par $U_{-\alpha}$.

Dans \cite{D1}, Demazure a associé à la décomposition $w=s_1\dots s_r$ de $w$ une résolution
des singularités $\pi:BS\fl \overline{X}_w$ (appelée résolution de Bott-Samelson),
où $BS=P_{\alpha_1}\times_B\dots\times_B P_{\alpha_r}/B$ et $\pi$ envoie la classe
$[p_1,\dots,p_r]$ de $(p_1,\dots,p_r)\in P_{\alpha_1}\times\dots\times P_{\alpha_r}$ dans
$BS$ sur la classe de $p_1\dots p_r$ dans $\overline{X}_w\simeq\overline{BwB}/B$.
La restriction de $\pi$ à $\pi^{-1}(X_w)$ induit un isomorphisme $\pi^{-1}(X_w)\fl X_w$.

On fait agir $B$ sur $BS$ par multiplication à gauche sur le premier facteur. (Le morphisme
$\pi$ est alors $B$-équivariant.) En particulier, on a une action de $\Gr_m$ sur $BS$ (via
le cocaractère $\lambda:\Gr_m\fl T$), dont les points fixes sont les
$[\gamma_1,\dots,\gamma_r]$, $(\gamma_1,\dots,\gamma_r)\in\Gamma$. Dans \cite{H}, Härterich a 
décrit la décomposition en cellules contractantes de $BS$. Nous allons rappeler ses
résultats.

Pour tout $(\gamma_1,\dots,\gamma_r)\in\Gamma$, le morphisme
$u_\gamma=(p_{\gamma_1(-\alpha_1)},\dots,p_{\gamma_r(-\alpha_r)}):\Aff^r\fl BS$
est une immersion ouverte (cf \cite{H}, $\S$ 1). On note $U_\gamma$ son image. On a 
$U_{(s_1,\dots,s_r)}=\pi^{-1}(X_w)$ et, pour tout $\gamma=(\gamma_1,\dots,\gamma_r)\in\Gamma$,
$U_\gamma\cap U_{(s_1,\dots,s_r)}=u_\gamma(\{(x_1,\dots,x_r)\in\Aff^r,x_i\not=0\ si\ \gamma_i=1\}$.

Soit $\gamma=(\gamma_1,\dots,\gamma_r)\in\Gamma$. Härterich (\cite{H}, formule 1.3) 
montre que la cellule contractante associée à $[\gamma_1,\dots,\gamma_r]$ est
$$C_\gamma=u_\gamma(\{(x_1,\dots,x_r)\in\Aff^r,x_i=0\ si\ i\not\in J(\gamma)\}).$$
La même méthode permet de montrer que la cellule dilatante associée à $\gamma$ est
$$C^\gamma=u_\gamma(\{(x_1,\dots,x_r)\in\Aff^r,x_i=0\ si\ i\in J(\gamma)\}).$$

La variété $BS$ est union disjointe des sous-variétés localement fermées $C^\gamma$,
$\gamma\in\Gamma$, et le morphisme $\pi:BS\fl\overline{X}_w$ est un isomorphisme au-dessus
de $X_w$, donc on a
$$X_w\simeq\bigcup_{\gamma\in\Gamma}^{\scriptscriptstyle{\bullet}}
U_{(s_1,\dots,s_r)}\cap C^\gamma.$$
Pour tout $\gamma\in\Gamma$, on pose
$$Y_\gamma=U_{(s_1,\dots,s_r)}\cap C^\gamma.$$

Soit $\gamma\in\Gamma$. D'après les formules ci-dessus pour $C^\gamma$
et $U_{(s_1,\dots,s_r)}\cap U_\gamma$, on a
$$Y_\gamma=u_\gamma(\{(x_1,\dots,x_r)\in\Aff^r,x_i\not=0\ si\ i\not\in I(\gamma)\ et\ x_i=0\ si\ i\in J(\gamma)\}).$$
Donc $Y_\gamma=\varnothing$ si $J(\gamma)\not\subset I(\gamma)$ et, si
$J(\gamma)\subset I(\gamma)$, 
$$Y_\gamma\simeq\Aff^{card(I(\gamma)-J(\gamma))}\times\Gr_m^{r-card(I(\gamma))}.$$

Enfin, comme le morphisme $\pi$ est $T$-équivariant, on a
$$\pi^{-1}(\overline{X}_w\cap X^v)=\bigcup_{\gamma\in\Gamma_v}C^\gamma,$$
d'où la deuxième formule de la proposition.

\section{Calcul des polynômes de Kazhdan-Lusztig en fonction des polynômes $R$}

Pour tout $d\in\Q$,
on définit un endomorphisme $\Q$-linéaire $\tau_{\leq d}$ de l'anneau des
polynômes de Laurent $\Q[t^{1/2},t^{-1/2}]$ par :
$$\tau_{\leq d}\left(\sum_{i\in\Z}a_it^{i/2}\right)=\sum_{i\leq d}a_it^{i/2}.$$

Le but de cette section est de donner une preuve géométrique du résultat
suivant, qui a été
prouvé de manière combinatoire par Brenti (cf \cite{B}, théorème 4.1) :

\begin{theoreme}\label{calcul_P_R} Avec les notations de la section 1, 
pour tous $v,w\in W$ tels que 
$v\leq w$, le polynôme $P_{v,w}$ est égal à :
\begin{flushleft}$\displaystyle{\tau_{\leq\ell(w)-\ell(v)-1}\left(\sum_{v=v_1<\dots<v_r<w}
(-1)^r R_{v_1,v_2}\tau_{\leq\ell(w)-\ell(v_2)}(R_{v_2,v_3}\dots\right.}$\end{flushleft}
\begin{flushright}$\displaystyle{\left.\dots \tau_{\leq\ell(w)-\ell(v_{r-1})}(R_{v_{r-1},v_r}\tau_{\leq\ell(w)-\ell(v_r)}(R_{v_r,w}))\dots)\right).}$\end{flushright}

\end{theoreme}

Pour tout schéma lisse $X$ sur $\Fi_q$, notons $\Tate(X)$ la sous-catégorie
triangulée de $D^b_m(X,\Q_\ell)$ engendrée par les objets isomorphes aux faisceaux
$\Q_\ell(m)$, $m\in\Z$. Le groupe de Grothendieck $K(\Tate(X))$ de $\Tate(X)$ est le groupe 
abélien libre engendré par les $[\Q_\ell(m)]$, $m\in\Z$. On définit un isomorphisme de groupes
$\varphi:\Tate(X)\fl \Z[t,t^{-1}]$ par $\varphi([\Q_\ell(m)])=t^{-m}$, pour tout $m\in\Z$.
Si $K$ est un objet de $\Tate(X)$ et $j\in\Nat^*$, on a, pour tout
$x\in X(\Fi_{q^j})$,
$$\varphi([K])(q^j)=Tr(F^{j*},K_{x}).$$

\begin{lemme}\label{proprietes_T}
\begin{itemize}
\item[(i)] Soit $X$ lisse connexe sur $\Fi_q$. Pour tout objet $K$ de $\Tate(X)$ et tout 
$a\in\Z$, on a :
$$\varphi([w_{\leq a}K])=\tau_{\leq a-dim(X)}(\varphi([K])).$$
\item[(ii)] Pour tous $v,w\in W$ tels que $v\leq w$ et tout $K\in\Tate(X_w)$, le complexe
$i_{v,w}^!j_{w!}K$ est dans $\Tate(X_v)$, 
on a un isomorphisme $Gal(\overline{\Fi}_q/\Fi_q)$-équivariant
$$(i_{v,w}^!j_{w!}K)_{x_v}\simeq R\Gamma_c((X_w\cap X^v)_{\overline{\Fi}_q},
K_{|X_w\cap X_v})$$ 
et $\varphi([i_{v,w}^!j_{w!}K])=\varphi([K])R_{v,w}(t)$.
\item[(iii)] Pour tous $v,w\in W$ tels que $v\leq w$, le complexe 
$i_{v,w}^*IC_{\overline{X}_w}$ est dans $\Tate(X_v)$.

\end{itemize}
\end{lemme}

\begin{preuve}\begin{itemize}
\item[(i)] On peut supposer que $K=\Q_\ell(m)$. L'égalité est alors évidente.
\item[(ii)] On peut supposer que $K=\Q_\ell(m)$. Le complexe $i_{v,w}^!j_{w!}K$ est un
complexe $B$-équivariant sur $X_v$. Comme $B$ agit transitivement sur $X_v$ et que les
stabilisateurs des points de $X_v$ sont connexes, les faisceaux de cohomologie de
$i_{v,w}^!j_{w!}K$ sont constants. Il suffit donc de montrer que
$(i_{v,w}^!j_{w!}K)_{x_v}\in\Tate(x_v)$.

D'après \cite{KL2} 1.4, on a un diagramme commutatif $T$-équivariant
$$\xymatrix@C=35pt{ & \overline{X_w} & X_w\ar[l]_-{j_w} \\
X_v\ar[ru]^-{i_{v,w}}\ar[r]_-{(id,x_v)} & X_v\times (\overline{X}_w\cap X^v)\ar[u] & X_v\times (X_w\cap X^v)\ar[u]\ar[l]^-{(id,j)}}$$
où les flèches verticales sont des immersions ouvertes et $j$ est l'inclusion
$X_w\cap X^v\fl \overline{X}_w\cap X^v$. D'après la formule de Künneth (cf SGA 5 III 1.7),
on a
$$i_{v,w}^!j_{w!}K=(id,x_v)^!(id,j)_!(\Q_{\ell,X_v}\boxtimes\Q_{\ell,X_w\cap X^v}(m))
=\Q_{\ell,X_v}\boxtimes (x_v^!j_!\Q_\ell(m)),$$
donc
$$(i_{v,w}^!j_{w!}K)_{x_v}=x_v^!j_!\Q_\ell(m)=x_v^!j_{w!}(K_{|X_w\cap X^v}).$$
D'après le sous-lemme \ref{lemme_homotopie} 
(qu'on peut appliquer par \cite{KL2} 1.5), on a un isomorphisme
$Gal(\overline{\Fi}_q/\Fi_q)$-équivariant
$$(i_{v,w}^!j_{w!}K)_{x_v}=x_v^!j_!(K_{|X_w\cap X^v})
\simeq R\Gamma_c((X_w\cap X^v)_{\overline{\Fi}_q},K_{|X_w\cap X^v}).$$
De plus, $X_w\cap X^v$ est union disjointe de schémas de la forme $\Aff^r\times\Gr_m^s$
d'après la proposition \ref{structure_link}, donc on a
bien $(i_{v,w}^!j_{w!}K)_{x_v}\in\Tate(x_v)$.

Enfin, l'isomorphisme ci-dessus donne, pour tout $j\in\Nat^*$,
$$\varphi([i_{v,w}^!j_{w!}K])(q^j)=Tr(F^{j*},(i_{v,w}^!j_{w!}K)_{x_v})
=Tr(F^{j*},R\Gamma_c((X_w\cap X^v)_{\overline{\Fi}_q},\Q_\ell(m))).$$
D'après la formule des points fixes de Grothendieck-Lefschetz (SGA 4 1/2 Rapport 1.3.2)
et le sous-lemme \ref{interpretation_R_1},
$$\begin{array}{rcl}\displaystyle{Tr(F^{j*},R\Gamma_c((X_w\cap X^v)_{\overline{\Fi}_q},
\Q_\ell(m)))} & = & \displaystyle{\sum_{x\in (X_w\cap X^v)(\Fi_{q^j})}Tr(F_x^*,\Q_\ell(m))} \\
 & = & \displaystyle{q^{-jm}card((X_w\cap X^v)(\Fi_{q^j}))=q^{-jm}R_{v,w}(q^j).}\end{array}$$
Donc, pour tout $j\in\Nat^*$,
$$\varphi([i_{v,w}^!j_{w!}K])(q^j)=q^{-jm}R_{v,w}(q^j)=\varphi([K])(q^j)R_{v,w}(q^j),$$
ce qui implique la dernière égalité de (ii).

\item[(iii)] Le complexe $i_{v,w}^*IC_{\overline{X}_w}$ est un complexe $B$-équivariant
sur $X_v$, donc ses faisceaux de cohomologie sont constants, et il suffit comme dans (ii)
de montrer que $IC_{\overline{X}_w,x_v}\in\Tate(x_v)$. Ceci résulte du théorème
\ref{calcul_P_IC}.

\end{itemize}
\end{preuve}

Le premier sous-lemme est prouvé par Springer dans \cite{S}, $\S$ 3, proposition 1.

\begin{souslemme}\label{lemme_homotopie} Soit $X$ une variété sur $\Fi_q$ munie d'une action de $\Gr_m$ qui contracte
$X$ sur un point $a\in X$ (c'est-à-dire que, pour tout $x\in X$, 
$\displaystyle{\lim_{\lambda\fl 0}\lambda.x=a}$). Soit $K$ un complexe $\ell$-adique
$\Gr_m$-équivariant sur $X$. Alors le morphisme d'adjonction $a_!a^!K\fl K$ induit un
isomorphisme
$$a^!K\iso R\Gamma_c(X,K).$$

\end{souslemme}

Le deuxième sous-lemme est une observation de Kazhdan et Lusztig (\cite{KL1} lemmes A3 et A4, voir
aussi \cite{KL2} 4.6).

\begin{souslemme}\label{interpretation_R_1} Pour tout $v\leq w$, pour tout $j\in\Nat^*$, on a
$$R_{v,w}(q^j)=card((X_w\cap X^v)(\Fi_{q^j})).$$

\end{souslemme}

\begin{preuvet}{\ref{calcul_P_R}}
Avec les notations ci-dessus, le théorème \ref{calcul_P_IC} implique :
pour tous $v,w\in W$ tels que $v\leq w$, on a 
$$P_{v,w}(t)=\varphi([i_{v,w}^*IC_{\overline{X}_w}]).$$
Le résultat cherché résulte alors directement du théorème \ref{calcul_IC} et du lemme
\ref{proprietes_T}.

\end{preuvet}

\end{document}